\documentclass[12pt,reqno]{amsart}

\usepackage[utf8]{inputenc}
\usepackage[normalem]{ulem}
\usepackage{amssymb}
\usepackage{amsmath}
\usepackage{amsthm}
\usepackage{amsfonts}
\usepackage{bbm}
\usepackage{bookmark}
\usepackage{hyperref}
\hypersetup{pdfstartview={FitH}}
\usepackage{mathrsfs} 
\usepackage{caption}

\theoremstyle{plain}
\newtheorem{theorem}{Theorem}[section]
\newtheorem{lemma}[theorem]{Lemma}

\theoremstyle{definition}

\theoremstyle{remark}

\numberwithin{equation}{section}
\newtheorem*{theorem*}{Theorem} 

\newcommand{\Z}{{\mathbb Z}}
\newcommand{\R}{{\mathbb R}}
\newcommand{\N}{{\mathbb N}}
\newcommand{\C}{{\mathbb C}}

\DeclareFontFamily{U}{mathx}{\hyphenchar\font45}
\DeclareFontShape{U}{mathx}{m}{n}{
<5> <6> <7> <8> <9> <10>
<10.95> <12> <14.4> <17.28> <20.74> <24.88>
mathx10
}{}
\DeclareSymbolFont{mathx}{U}{mathx}{m}{n}
\DeclareFontSubstitution{U}{mathx}{m}{n}
\DeclareMathAccent{\widecheck}{0}{mathx}{"71}

\begin{document}

\title[Norm-variation of cubic ergodic averages]{Norm-variation of cubic ergodic averages}

\author[P. Durcik]{Polona Durcik}
\address{Schmid College of Science and Technology,  Chapman University, One University Drive, Orange, CA 92866, USA}
\email{durcik@chapman.edu}

\author[K. A. \v{S}kreb]{Kristina Ana \v{S}kreb}
\address{Faculty of Civil Engineering, University of Zagreb, Fra Andrije Ka\v{c}i\'{c}a Mio\v{s}i\'{c}a 26, 10000 Zagreb, Croatia}
\email{kristina.ana.skreb@grad.unizg.hr}

\date{\today}


\begin{abstract}
We prove a quantitative result on norm convergence of   cubic ergodic averages with respect to $d\geq 1
$  commuting measure-preserving transformations. We use  harmonic analysis techniques, a key tool being   estimates for   singular Brascamp-Lieb forms with cubical structure, which are used as a black box. 
\end{abstract}

\maketitle

\section{Introduction}
Let   $(X,\mathcal{F},\mu)$  be   a probability space and $d\geq 1$. Let   $T_1,T_2,\dots,T_d:X\to X$ be  mutually  commuting, measure-preserving transformations; that is,   $T_iT_j=T_jT_i$ for every $i \neq j$, and $\mu(T_i^{-1}(E))=\mu(E)$ for every $i$ and $E\in\mathcal{F}$.
Let \[Q = \{0,1\}^d \setminus\{(0,\ldots, 0)\}.\]
For $n\geq 1$, a tuple $\mathbf{f}= (f_j)_{j\in Q}$ of    functions $f_j\in L^\infty(X)$, and $x\in X$,    define
\begin{equation}
M_n(\mathbf{f})(x)=\frac{1}{n^d}\sum_{i_1,\dots, i_d=0}^{n-1}\prod_{j\in Q} f_{j}\big(\big( \prod_{l=1}^dT_l^{j_li_l} \big) x\big),
\label{cub_av}
\end{equation}
where $j=(j_1,\ldots, j_d)$. 
We prove the following   norm-variation estimates. 
\begin{theorem}
\label{thm:ergodic}
Let   $d\ge 1$ and  $1\le p < \infty$. If  $d=1$,  let  $\varrho \ge \max\{2,p\}$ and if $d\geq 2$,  let 
\[\varrho>\max\{2,p(2^d-1)/2^{d-1}\}.\] 
There exists a constant $C>0$ such that for any  probability space $(X,\mathcal{F},\mu)$,  any    mutually commuting measure-preserving transformations  $T_1,T_2,\dots,T_d:X\to X$,  any positive integers  $I$ and $n_0< n_1< \cdots  <n_I$, and  any tuple of functions $\mathbf{f}=(f_j)_{j\in Q}$   with  $f_j \in  {L}^\infty(X)$ and $\|f_j\|_{L^\infty(X)} =1$ for  each   $j\in Q$,    
\begin{equation*}
\sum_{i=1}^I \|M_{n_i}(\mathbf{f})-M_{n_{i-1}}(\mathbf{f})\|^\varrho_{L^p(X)} \le  C.
\end{equation*}
\end{theorem}
The averages \eqref{cub_av} are also called cubic ergodic averages.  
Theorem \ref{thm:ergodic} establishes a quantitative form of norm convergence of these averages.  
Their qualitative $L^2$ convergence was proved in \cite{Au10} and \cite{H09} using different methods. In the special case when $T_1=T_2=\cdots=T_d$, the pointwise convergence of \eqref{cub_av} for a.e. point $x\in X$ was proved in \cite{A10}, \cite{CF12}, and \cite{HSY14}. The pointwise convergence in the case of general commuting transformations was established for two transformations in \cite{DS16} and extended to more transformations   in \cite{DS18}.

Cubic ergodic averages \eqref{cub_av} appeared for instance in \cite{H09} and in the proof of the $L^2$ convergence of multiple ergodic averages 
\begin{equation}
\frac{1}{n}\sum_{i=0}^{n-1}f_1(T_1^ix)f_2(T_2^ix)\cdots f_d(T_d^ix).
\label{mea}
\end{equation}
The case $d=1$ of \eqref{mea} is classical, while their study for $d\ge 2$
was motivated by the works \cite{F77}, \cite{FK78}, and \cite{FKO82}, and it influenced the development of various tools in ergodic theory and topological dynamics. The norm convergence of \eqref{mea} was shown for all $d\geq 1$ in \cite{T08}, 
reproved in \cite{Au10}, \cite{H09} and once again in \cite{W12} in the more general case when the transformations generate a nilpotent group. The case $d=2$ has been known previously \cite{CL84}. 
Almost everywhere convergence of \eqref{mea} is still an open problem when $d\ge 2$. For some partial progress on that matter we refer the interested reader to \cite{DS16} and \cite{DS18}. 

Using techniques from harmonic analysis, quantitative $L^2$ convergence of \eqref{mea} was established for $d=2$ and $d=3$ in \cite{DKST16} and \cite{DST23}, respectively. 
These papers prove norm variation estimates for \eqref{mea}, sharp in the case $d=2$. The result for $d=2$ was previously shown in \cite{K16} in a  simplified model.
The papers \cite{DKST16} and \cite{DST22} rely on harmonic analysis methods, where the key ingredients are certain $L^p$ and   cancellation estimates  for operators related to   the triangular Hilbert transform. The analogues of these estimates have not yet been established for  $d\geq 4$. Any quantitative results for norm convergence of \eqref{mea} for $d\geq 4$ remain an open problem.

In contrast to this, for the cubic averages \eqref{cub_av}  we are able to prove norm-variation estimates for all $d\ge 1$. We  also approach the problem via real harmonic analysis, however, here estimates for less singular forms suffice. Namely, the main ingredients in our proof   are $L^{2^d}$-estimates for multilinear singular integrals with cubical structure from \cite{DST22}, which are known in all dimensions. The key restriction on the variation exponent $\varrho$ in Theorem \ref{thm:ergodic} stems from these estimates. It would be interesting to study sharp norm-variation bounds for the averages \eqref{cub_av}. 
We note that the harmonic analysis techniques used here are not sufficient to reprove almost everywhere convergence of the cubic ergodic averages. Quantifying their pointwise convergence remains an interesting open problem.

To prove Theorem \ref{thm:ergodic}, we   introduce the following   variant of the averages on the Euclidean space and prove   estimates for them. For a tuple $\mathbf{f}=(f_j)_{j\in Q}$ of functions $f_j \in  {L}^{2^d}(\R^d)$   and $t>0$, we define 
\[ A_t (\mathbf{f})(x) = t^{-d}\int_{[0,t)^d} \prod_{j\in Q}f_{j}(x+j\odot s)\,ds,\]
where $j\odot s=(j_1s_1,\dots,j_ds_d)$.
We denote the $L^p$ norm of a function $f$ on $\R^d$ by
\[\|f\|_p = \|f\|_{L^p(\R^d)}.\]
Throughout the paper, the letter $q$ will always denote the exponent
\[q=2^d/(2^d-1).\]
Its     H\"older conjugate is $q'=2^{d}$, i.e. $q^{-1} + {(q')^{-1}} =1$.
We prove the following theorem. 
\begin{theorem}\label{thm:char_q}
For every $d\ge 1$  
there exists a constant $C>0$ such that the following holds. For any positive integer $I$, any positive real numbers $t_0< t_1< \cdots  <t_I$, and any tuple $\mathbf{f}=(f_j)_{j\in Q}$ of functions $f_j \in  {L}^{2^d}(\R^d)$  with $\|f_{j}\|_{2^d}=1$ for each $j\in Q$,  
\begin{equation}\label{char_var}
     \sum_{i=1}^{I} \big\|A_{t_i}(\mathbf{f})-A_{t_{i-1}}(\mathbf{f})\big\|_{L^{q}(\R^d)}^{q}   
\le C I^{1-\frac{q}{2}} .
\end{equation}
\end{theorem}

The proof of Theorem \ref{thm:char_q} is presented in Section \ref{sec:proofmain}. The idea is to split the norm-variation  into the {long variation}, corresponding to the  dyadic scales, and the {short variation}, corresponding to the scales from a fixed dyadic interval. This is the usual approach to variational estimates, as can be seen, for instance, in \cite{JSW08} and \cite{DOP17}. The long and short variation estimates are then shown by   a lacunary decomposition of the characteristic function of the unit cube and carefully estimating the resulting pieces. The factor $I^{1-q/2}$ arises when passing to a square function,  whose bounds are  in turn reduced to  a multilinear singular integral estimate from \cite{DST22}.

To pass from Theorem \ref{thm:char_q} to Theorem \ref{thm:ergodic}, in Section \ref{sec:transitionergodic} we  perform a number of   standard reductions. 
The  estimate in Theorem \ref{thm:char_q}  implies that the sequence of averages has at most $O(\varepsilon^{-2})$ jumps of size $\varepsilon$ in the $L^q$ norm. 
A layer-cake decomposition of the $\ell^\varrho$ norm then gives $\varrho$-variation estimates in $L^q$ for $\varrho>2$.  The  transition to ergodic averages from Theorem \ref{thm:char_q}  ultimately follows by a variant of Calder{\'o}n's transference principle. To pass to the other exponents $p$ we  use monotonicity and log-convexity of $L^p$ norms. \\

{\bf Acknowledgment}
 P.~D.  was partially supported by the grant  NSF DMS-2154356 and by a grant from
the Simons Foundation SFI-MPS-TSM-00013943.  K.~A.~\v{S}. was partially supported by the Croatian Science Foundation under the project number HRZZ-IP-2022-10-5116 (FANAP). We thank Vjekoslav Kova\v c   and Christoph Thiele for inspiring discussions.  The second author thanks Pavel Zorin-Kranich for his comments on an earlier version of the paper.

\section{Proof of Theorem \ref{thm:char_q} using singular integral estimates}
\label{sec:proofmain}
Throughout this and the subsequent sections, the symbol $C$ will denote various positive constants, possibly changing from line to line. Unless stated otherwise, such constants may depend on the dimension $d\ge 1$ and on the Schwartz seminorms of the bump-like functions appearing in the lemmas below. For a function
$\rho:\R^n\to \C$, $n\ge 1$, we will use the notation
\[\rho_{(t)}(u)=t^{-n}\rho(t^{-1}u).\]

 The first step in the proof of Theorem \ref{thm:char_q} is  to  decompose the characteristic function of the cube. We use a similar one-dimensional decomposition as in \cite{DKST16}, applied to   each side of the cube.  
 
Let $\chi$ and $\vartheta$ be  Schwartz functions such that  $\widehat{\chi}$ is  even, non-negative, smooth,           supported in $[-1,1]$,  identically equal to one on $[-1/2,1/2]$, and such that 
\[(\widehat{\vartheta})^2=\widehat{\chi}-\widehat{\chi_{(2)}}.\]
We denote \[\theta = \vartheta * \vartheta.\]
Then $\widehat{\vartheta}$ and $\widehat{\theta}$  are   supported in $[-1,-1/4]\cup [1/4,1]$, $\vartheta$ is real-valued,  
and for all $\xi\neq 0$,
\[\sum_{k\in \Z}\widehat{\theta}(2^k\xi) =1.\]
By $\widetilde{\theta}$ and  $\widetilde{\vartheta}$   we denote the primitives of $\theta$ and ${\vartheta}$, respectively. Then $\widehat{\widetilde{\vartheta}}$, and $\widehat{\widetilde{\theta}}$ are also  supported in $[-1,-1/4]\cup [1/4,1]$  and
\[\widetilde{\theta}(x) = \int_{-\infty}^x\theta(s)ds =  (\widetilde{\vartheta}*\vartheta)(x). \]
We denote
\[\phi = \mathbf{1}_{[0,1)}\ast\chi\]
and for $m\le -1$ also 
\[\phi_{0,m}= \widetilde{\theta}(2^{-m}\cdot ),\]
\begin{equation}
\label{phi1m}
    \phi_{1,m} =  -\widetilde{\theta}(2^{-m}(\cdot-1)),
\end{equation}
The convolution structures of $\theta$ and $\widetilde{\theta}$ will be used   in some of the arguments  below,  when analyzing the averages associated with $\phi_{1,m}$.

We decompose
\begin{equation*}
   \mathbf{1}_{[0,1)} = \mathbf{1}_{[0,1)}\ast\chi + 
\sum_{m=-\infty}^{-1} \mathbf{1}_{[0,\infty)}\ast\theta_{(2^m)}
-\sum_{m=-\infty}^{-1} \mathbf{1}_{[1,\infty)}\ast\theta_{(2^m)} 
\end{equation*}
\[=\phi + \sum_{m=-\infty}^{-1}\phi_{0,m} + \sum_{m=-\infty}^{-1} \phi_{1,m}.\]
This decomposition holds in a pointwise a.e. sense, which follows, for instance, from the weak $L^2$ boundedness of the maximally
truncated convolution-type singular integrals.  

Applying this decomposition to each factor of $\mathbf{1}_{[0,1)^d}$, we obtain
\begin{equation}
    \label{char_decomp}
    \mathbf{1}_{[0,1)^d}(s) = \prod_{l=1}^d\mathbf{1}_{[0,1)}(s_l) = \sum_{n=0}^d \sum_{\substack{S \subset \{1,\ldots,d\} \\ |S| = n}}  \sum_{\epsilon_1,\ldots, \epsilon_n\in\{0,1\}}\sum_{m_1,\ldots,m_n=-\infty}^{-1} \Phi_{n,S,\epsilon,m} (s),
\end{equation}
where $m=(m_1,\ldots, m_n),\epsilon=(\epsilon_1,\ldots, \epsilon_n)$, $S=\{l_1<\cdots < l_n\}$, and write 
\[\Phi_{n,S,\epsilon,m}(s) = \prod_{l=1}^d \varphi_{l,n,S,\epsilon,m}(s_l),\]
where $\varphi_{l,n,S,\epsilon,m} = \phi_{\epsilon_i,m_i}$ if $l=l_i$ for some $1\le i\le n$, and $\varphi_{l,n,S,\epsilon,m} =\phi$ otherwise.
In words, $n$ is the number of coordinates where we do not use the function $\phi$, and the set $S$ gives the specific set of these coordinates. The numbers $\epsilon_i$ indicate whether we pick $\phi_{0,m_i}$ or $\phi_{1,m_i}$, where $m_i$ is the dyadic scale of that function.

For an integrable function $\Phi$ on $\R^d$, a tuple $\mathbf{f}=(f_j)_{j\in Q}$ of functions $f_j \in  {L}^{2^d}(\R^d)$, and $t>0$, we define 
\begin{equation*}
A_t^\Phi (\mathbf{f})(x)=\int_{\R^d} \Big( \prod_{j\in Q}f_{j}(x+j\odot s) \Big) \Phi_{(t)}(s)\,ds.
\end{equation*}
  When $\Phi = \mathbf{1}_{[0,1)^d}$, this definition recovers the averages $A_t (\mathbf{f})$.

For each $n,S,\epsilon$ in \eqref{char_decomp} we will show that   there exists  $C>0$ such that for any tuple of the negative integers  $m=(m_1,\ldots, m_n)$,   positive integer $I$,   sequence $(t_i)_{i=0}^I$ of  positive real numbers, and   tuple  $\mathbf{f}$ as in \eqref{char_var},  the following holds: if $\Phi  =  \Phi_{n,S,\epsilon,m}$, then  
\begin{equation}
    \label{est_decay}
     \sum_{i=1}^{I} \big\|A_{t_i}^{\Phi}(\mathbf{f})-A_{t_{i-1}}^{\Phi}(\mathbf{f})\big\|_{q}^{q} 
\le C 2^{\delta  m'}  I^{1-\frac{q}{2}},
\end{equation}  
where $\delta =0$ if $n=0$ and $\delta=q/q'$ if $n>0$, and $m'=\min(m_1,\ldots,m_n)$. 
Using the decomposition \eqref{char_decomp}, applying the triangle inequality for  the $L^q$ and $\ell^q$ norms, taking the largest constant $C$ over $n,S,\epsilon$,  and summing in all parameters $n,S,\epsilon,m$  we then obtain  Theorem~\ref{thm:char_q}.

By  a change of variables, permutation of  the functions $f_j$, and permutation of the arguments in  $f_j$, 
we may assume that $\Phi$ is arranged so that its factors appear in the following order: first $\phi$, then $\phi_{0,m}$, and finally $\phi_{1,m}$.  
That is, it suffices to prove that for any integers $0\le L_1\le L_2\le d$, there is $C>0$ such that for any   tuple of the negative integers $(m_l)_{l=L_1+1}^d$, and any 
$I$,     $(t_i)_{i=0}^I$,   and     $\mathbf{f}$ as in \eqref{char_var}, the following holds. Let      
\begin{equation}
\label{phi_gen}
    \Phi(s) = \prod_{l=1}^d\varphi_{l,m_l}(s_l),
\end{equation}
where 
 $\varphi_{l,m_l} =\phi$ for $1\le l \le L_1$, $\varphi_{l,m_l} = \phi_{0,m_l}$ for $L_1+1\le l\le L_2$,  and $\varphi_{l,m_l} = \phi_{1,m_l}$ for $L_2+1\le l \le d$. If $L_1=d$, let  $\delta=0$  and if $L_1<d$, let  $\delta=q/q'$, and let $m' = \min(m_{L_1+1},\ldots,m_{d})$. 
Then it suffices to show that \eqref{est_decay} holds.
 
 Note that in the case $L_1<d$    there is at least one mean-zero function $\phi_{0,m_l}$ or $\phi_{1,m_l}$, while in the case $L_1=d$,  all functions are equal to $\phi$.

Using a standard separation into long and short jumps (see e.g. \cite{JSW08}), it suffices to establish the long and short variation bounds 
\begin{equation}
    \label{eq:long}
    \sum_{i=1}^I \|A^\Phi_{2^{k_i}}(\mathbf{f})-A^\Phi_{2^{k_{i-1}}}(\mathbf{f})\|_q^q\leq C 2^{\delta m'} I^{1-\frac{q}{2}} ,
\end{equation}
\begin{equation}
    \label{eq:short}
    \sum_{i=0}^I \sup_{\substack{n\in \Z_+\\ 2^{k_i} < t_0<\ldots<t_n\le 2^{k_{i}+1}}}  \sum_{\ell=1}^n \|A^\Phi_{t_\ell}(\mathbf{f})-A^\Phi_{t_{\ell-1}}(\mathbf{f})\|_q^q\leq C  2^{\delta m'} I^{1-\frac{q}{2}}
\end{equation}
for any   increasing sequence of the integers $(k_i)_{i=0}^{I}$, with $C$ independent of $m,I,\mathbf{f},(k_i)_{i=0}^I$.

Before proceeding with the proofs we  formulate a number of preparatory results.
Throughout the rest of the paper we denote
\[Q_0=\{0,1\}^d.\]
For a tuple $\mathbf{F}=(F_j)_{j\in Q_0}$ of Schwartz functions $F_j:\R^d\to \C$ and a tempered distribution $K$ we define the singular integral form
\[\Lambda(K,\mathbf{F}) =  \textup{p.v.} \int_{\R^{2d}} \Big( \prod_{j\in Q_0} F_j(x+j\odot s) \Big)   K(s) \, ds\, dx.\] 
The results below  will  rely on   the following estimate for singular Brascamp-Lieb forms with cubical structure, which is  a special case of Theorem 1.1 in \cite{DST22}. 
\begin{theorem}[Theorem 1.1 in \cite{DST22}]
\label{thm:singint} Let $d\ge 1$. There is a constant $C>0$ such that the following holds. 
Let $K$ be a tempered distribution  such that $\widehat{K}$ is a smooth function on $\R^d\setminus\{0\}$ and satisfies 
\begin{equation}
    \label{symest}
    |\partial^\alpha \widehat{K}(\xi)|\le |\xi|^{-|\alpha|}
\end{equation}
for all multi-indices $|\alpha|\le 2^{6d}$ and all $\xi\neq 0$.
Then for all tuples $\mathbf{F}=(F_j)_{j\in Q_0}$ of Schwartz functions $F_j:\R^d\to \C$ with  $\|F_j\|_{2^d} =1$ for all $j\in Q_0$, 
    \begin{equation*}
                    | \Lambda(K,\mathbf{F})  | \le C.
    \end{equation*}
\end{theorem}
While this result is stated only for 
$d\geq 2$ in \cite{DST22},  the $d=1$ case reduces to the classical convolution setting:
here $\Lambda$ is a bilinear form associated with a classical convolution-type operator,  and the bound on $L^2$  follows by the Cauchy-Schwarz inequality and Plancherel's theorem.

 We will use Theorem \ref{thm:singint} for kernels of the form
 \begin{equation}
     \label{K_apply}
     K =  \sum_{k=k_0}^{k_1}\varepsilon_k \Psi_{(2^k)}
 \end{equation}
where $k_0,k_1$ are integers, $|\varepsilon_k|\le 1$,  $\Psi$ is smooth, $\widehat{\Psi}\in C^{100^d}$,   $\widehat{\Psi}(0)=0$. Such kernels satisfy  \eqref{symest} up to a multiplicative constant independent of $k_0,k_1$. To see that,  we differentiate and use the triangle inequality, giving  
\[|\partial^{\alpha}\widehat{K}(\xi)| = \Big|\sum_{k=k_0}^{k_1}\varepsilon_k \partial^{\alpha}\widehat{\Psi}(2^k\xi)\Big|  \le  \sum_{k=k_0}^{k_1}2^{k|\alpha|} |(\partial^\alpha \widehat{\Psi})(2^k\xi)| \]
Let $l\in \Z$ be such that $2^l\le |\xi|\le 2^{l+1}$. If $k\le -l$, we use that $\widehat{\Psi}$   vanishes at the origin,  and thus $|(\partial^\alpha \widehat{\Psi})(2^k\xi)| \le C 2^{k}|\xi|\le C 2^{k+l}$. This yields   
\[  2^{k|\alpha|} |(\partial^\alpha \widehat{\Psi})(2^k\xi)| \le   C2^{-l|\alpha|} 2^{k+l} \le C |\xi|^{-|\alpha|}2^{k+l}. \]
If $k\ge -l$, we estimate  using the rapid decay of $\widehat{\Psi}$,
\[ 2^{k|\alpha|} |(\partial^\alpha \widehat{\Psi})(2^k\xi)| \le C 2^{k|\alpha|}  (2^k|\xi|)^{-|\alpha|-1}   \le C |\xi|^{-|\alpha|} 2^{-(k+l)}.\]
  This gives for any $|\alpha|\le 2^{6d}$
  \[   |\partial^{\alpha}\widehat{K}(\xi)|  \le C \Big( |\xi|^{-|\alpha|}\sum_{k\le -l}2^{k+l} + |\xi|^{-|\alpha|}\sum_{k\geq -l}2^{-(k+l)} \Big)\le C_{\Psi,d}|\xi|^{-|\alpha|},\]
  where $C_{\Psi,d}$ is  a constant that depends on $\Psi$ and $d$. Theorem \ref{thm:singint} then applies to the kernel $C_{\Psi,d}^{-1}K$.

The estimate in Theorem \ref{thm:singint} is invariant under anisotropic rescalings of the form $\Lambda$, which will be used multiple times to simplify the arguments below. 
\begin{lemma}
\label{lem:rescalingtp}
Let $K\in L^1(\R^d)$. Assume that for any tuple   $\mathbf{F}=(F_j)_{j\in Q_0}$  of Schwartz functions $F_j:\R^d\to \C$ with  $\|F_j\|_{2^d} =1$ for all $j\in Q_0$,  
    \begin{equation}
        \label{ani_scal_aspt}
        | \Lambda(K,\mathbf{F})  | \le  1.
    \end{equation}
    For $a=(a_1,\ldots, a_d),\, a_i\in \R$, let $K_a(s) = a_1^{-1}\ldots a_d^{-1} K(a_1^{-1}s_1,\ldots, a_d^{-1}s_d)$. Then for any tuple   $\mathbf{F}=(F_j)_{j\in Q_0}$  of Schwartz functions $F_j:\R^d\to \C$ with  $\|F_j\|_{2^d} =1$ for all $j\in Q_0$,  
    \begin{equation*}
         | \Lambda(K_a,\mathbf{F})  | \le  1.
    \end{equation*}
\end{lemma}

\begin{proof}  By the change of variables  $x\to a\odot x$ and $s \to a \odot s$,  
    \[ \Lambda(K_a,\mathbf{F}) =   a_1\cdots a_d\int_{\R^{2d}} \Big( \prod_{j\in Q_0} F_j(a\odot x+j\odot a \odot  s ) \Big)   K(s) dx ds  \]
    \[=  \int_{\R^{2d}} \Big( \prod_{j\in Q_0} F_{j,a}(x+j\odot s) \Big)   K(s) dx ds, \]
where $F_{j,a}(y) = (a_1\cdots a_d)^{2^{-d}} F_j(a\odot y) $. For each $j\in Q_0$,  \[\|F_{j,a}\|_{2^d} = \|F_j\|_{2^d} = 1.\]  
Applying \eqref{ani_scal_aspt} to the tuple $(F_{j,a})_{j\in Q_0}$  yields the claim.
\end{proof}

The following result is an application of Khintchine's inequality. 
\begin{lemma} \label{lem:kintchine}  
There exists a constant $C>0$ such that the following holds for every positive integer $I$ and every tuple  $(\Phi_i)_{i=1}^I$   of integrable functions $\Phi_i:\R^d\to \R$. If  for any     real numbers  $\varepsilon_i$  with   $|\varepsilon_i|\le 1$    and any  tuple 
 $\mathbf{F}=(F_j)_{j\in  Q_0}$ of Schwartz functions $F_j:\R^d\to \R$ with $\|F_j\|_{2^d}=1$ one has
\begin{equation}
    \label{singint_assumption}
    \Big |\Lambda \Big(\sum_{i=1}^I \varepsilon_i \Phi_{i},\mathbf{F} \Big ) \Big |\le 1,   
\end{equation}
then   for any tuple  $\mathbf{f}=(f_j)_{j\in  Q}$ of functions  $f_j\in L^{2^d}(\R^d)$ with $\|f_j\|_{2^d}=1$, 
\[\sum_{i=1}^I  \|A^{\Phi_i}(\mathbf{f})\|_q^{q} \le C I^{1-q/2}.\]
  
\end{lemma}

    \begin{proof} Let $I,\Phi_i,\mathbf{f}$  be given   and set 
    \[A_i=A^{{\Phi_i}}(\mathbf{f}).\]
    Recall that $q =2^d/(2^d-1)\in [1,2]$.
     Using   
the power mean inequality   
\begin{equation*} 
\Big(\frac{1}{I}\sum_{i=1}^I |A_i|^{q}\Big)^{\frac{1}{q}} \le \Big(\frac{1}{I}\sum_{i=1}^I |A_i|^2\Big)^{\frac 12}.
\end{equation*} 
We obtain 
\[\sum_{i =1}^I \|A_i\|_q^q  = \Big\|\Big( \sum_{i =1}^I |A_i|^q \Big)^{\frac{1}{q}}\Big\|_q^q
\le  I^{1-\frac{q}{2}} \Big\| \Big( \sum_{i =1}^I |A_i|^2 \Big)^{\frac{1}{2}}\Big\|_q^q. \]
Thus, it remains to show 
\[\Big\| \Big( \sum_{i =1}^I |A_i|^2 \Big)^{\frac{1}{2}}\Big\|_q^q \le C.\]

Let $\varepsilon = (\varepsilon_1,\ldots,\varepsilon_I)$ where $\varepsilon_i\in \{-1,1\}$ are i.i.d.~random  signs. 
  By  Khintchine's inequality,  there is a constant $C_q$ depending only on $q$ such that 
\[     \Big( \sum_{i=1}^I  |A_i |^2  \Big)^{1/2}   \le  C_q\Big( \mathbf{E}  \Big| \sum_{i=1}^I \varepsilon_i A_i\Big|^q\Big)^{1/q},\]
Therefore, 
\[    \Big \| \Big( \sum_{i=1}^I  |A_i |^2  \Big)^{1/2}\Big\|^q_q   \le  C_q^q \Big\| \Big( \mathbf{E}  \Big| \sum_{i=1}^I \varepsilon_i A_i\Big|^q\Big)^{1/q}\Big\|^q_q \]
\[= C_q^q   \int_{\R^d}\mathbf{E}  \Big| \sum_{i=1}^I \varepsilon_i A_i(x)\Big|^q dx  = C_q^q\,  \mathbf{E}  \int_{\R^d}\Big| \sum_{i=1}^I \varepsilon_i A_i(x)\Big|^q dx . \]
It thus remains to show 
    \[\Big \| \sum_{i=1}^{I} {\varepsilon}_i A^{{\Phi_i}}(\mathbf{f}) \Big \|_q \le C    \]
and then average  over $\varepsilon$, which in turn finishes the proof.  

    To see this estimate, we dualize the $L^q$ norm, after which it suffices to show that  for any $f_0\in L^{q'} = L^{2^{d}}$ with $\|f_0\|_{2^d}=1$, 
            \[    \Big| \int_{\R^d}f_0(x)\Big(\sum_{i=1}^{I} {\varepsilon}_i A^{{\Phi_i}} (\mathbf{f}) (x)\Big) \, dx \Big|  \le C.\]
By density of the Schwartz functions in $L^{2^d}$,  it suffices to prove this for Schwartz functions $f_j$. By splitting $f_j$ into real and imaginary parts, we may assume they are real-valued.
Expanding out the definition of $A^{{\Phi_i}}$, the form in the last display can be recognized as the form in \eqref{singint_assumption}.
The claim now follows from the assumption \eqref{singint_assumption}. 
    \end{proof}

\subsection{Long variation: estimate \texorpdfstring{\eqref{eq:long}}{}}
First we  prove the estimate when $L_1=d$ in \eqref{phi_gen}.

\begin{lemma} 
 \label{lem:long_var} 
 Let $\Phi = \prod_{l=1}^d \phi$.
There exists  a constant $C>0$ such that for  any  positive integer   $I$, increasing sequence of the integers $(k_i)_{i=0}^I$, and $\mathbf{f}$ as in \eqref{char_var},   
         \[\sum_{i =1}^I \|A^\Phi_{2^{k_i}}(\mathbf{f}) - A^\Phi_{2^{k_{i-1}}}(\mathbf{f})\|_q^q \le C I^{1-\frac{q}{2}}.\]
 \end{lemma}

 \begin{proof} 
We have
\[A^\Phi_{2^{k_i}}(\mathbf{f}) - A^\Phi_{2^{k_{i-1}}}(\mathbf{f}) = A^{\Phi_{(2^{k_i})}-\Phi_{(2^{k_{i-1}})}}(\mathbf{f}).\]
It suffices to show that there is a constant $C_0$   such that for any $I$, $(k_i)_{i=0}^I$, any  $|\varepsilon_i|\le 1$,   and any tuple $\mathbf{F}$ of normalized Schwartz functions, the estimate
\begin{equation}
     \label{longvar1_singint}
     |\Lambda(K,\mathbf{F})| \le C_0 
 \end{equation}
holds with
\[K = \sum_{i=1}^I \varepsilon_i (\Phi_{(2^{k_i})}-\Phi_{(2^{k_{i-1}})}). \]
To finish the proof of   Lemma \ref{lem:long_var} one then applies  Lemma \ref{lem:kintchine}   to the tuple of functions $(\max(C_0,1)^{-1}(\Phi_{(2^{k_i})}-\Phi_{(2^{k_{i-1}})}))_{i=1}^I$ and uses homogeneity of the $L^q$ norm. 

To see \eqref{longvar1_singint}, we  write
\[    \sum_{i=1}^I \varepsilon_i (\Phi_{(2^{k_i})}-\Phi_{(2^{k_{i-1}})})   =    \sum_{i=1}^I \varepsilon_i \sum_{k=k_{i-1}+1}^{k_i} \Phi_{(2^{k})}-\Phi_{(2^{k-1})}   
= \sum_{k=k_0}^{k_I} \tilde{\varepsilon}_k (\Phi_{(2^{k})}-\Phi_{(2^{k-1})}),  \]
where $\tilde{\varepsilon}_k = \varepsilon_i$ for $k_{i-1}+1\le k \le k_i$. Thus, it suffices to show \eqref{longvar1_singint} 
with 
\[K = \sum_{k=k_0}^{k_I} \tilde{\varepsilon}_k \Psi_{(2^k)},\]
where  $\Psi = \Phi-\Phi_{(2^{-1})}$. 
This kernel is of the form \eqref{K_apply}. Therefore, the desired bound \eqref{longvar1_singint} follows from  
 Theorem~\ref{thm:singint}.
 \end{proof}

The case $L_1<d$ will involve various translation parameters $m_l$ appearing in \eqref{phi_gen}. In some cases they will be removed with an application of the Cauchy-Schwarz inequality as in  the following lemma.  
\begin{lemma} \label{lem:cs}  
For     $1\le l\le d$ let $\varphi_l,\rho_l:\R\to \R$ be integrable functions 
satisfying
\[|\varphi_l| \le \rho_l. \]
Assume that there exist integrable functions  $\sigma_0,\sigma_1:\R\to \R$  so that  for some  $1\le l_0\le d$,  
\[\varphi_{l_0} = \sigma_0 * \sigma_1.\] 
Define the kernels $K$ and $K_b$, $b\in \{0,1\}$, by
\[K(s)=\prod_{l=1}^d\varphi_l(s_l),\quad K_b(s) = (\sigma_b^\star*\sigma_b)(s_{l_0})\prod_{l=1,l\neq l_0}^d\rho_l(s_l), \] 
where $\sigma_b^\star(u) = \sigma_b(-u)$. 
Then for all tuples 
$\mathbf{F}=(F_j)_{j\in Q_0}$ of Schwartz functions $F_j:\R^d\to \R$, 
\[ | \Lambda(K,\mathbf{F})  |\le   \Lambda(K_0,\mathbf{F})^{1/2}\Lambda(K_1,\mathbf{F})^{1/2} . \]
\end{lemma}
\begin{proof}
 Without loss of generality we may assume $l_0=d$.
Let   $x^0=(x_1^0,\ldots, x_d^0)\in \R^d, x^1 = (x_1^1,\ldots, x_d^1)\in \R^{d}$, $j=(j_1,\ldots, j_d)$. We introduce new variables 
$x^0=x$ and $x^1=x+s$, and write $x^j =(x_1^{j_0},\ldots, x_d^{j_d})$. This gives
\begin{equation}
    \label{form_equiv}
      \Lambda(K,\mathbf{F})  = \int_{\R^{2d}} \Big( \prod_{j\in Q_0} F_j(x^j) \Big) K(x^1-x^0)   dx^0 dx^1. 
\end{equation}
We write
\[\varphi_d(x_{d}^1-x_{d}^0) =\int_{\R} \sigma_0(p) \sigma_1(x_{d}^1-x_{d}^0-p) dp =  \int_{\R} \sigma_0(p-x_{d}^0) \sigma_1(x_{d}^1-p) dp .\] 
Denoting $\psi_0 = \sigma_0^\star$ and $\psi_1=\sigma_1$, 
  we can  write the form \eqref{form_equiv} as
\[ \int_{\R^{2d-1}} \Big( \prod_{b=0}^1  \int_{\R} \Big( \prod_{j\in Q_0, j_d=b} F_j(x^j) \Big) \psi_b(x_{d}^b-p) dx_d^b \Big)\] 
\[\times \Big( \prod_{l=1}^{d-1} \varphi_l(x_l^1-x_l^0)\Big) (dx^0_l)_{l=1}^{d-1}(dx^1_l)_{l=1}^{d-1} dp.\]
We use the triangle inequality and estimate $|\varphi_l|\le \rho_l$ for $1\le l \le d-1$ to bound the absolute value of this display by 
\[ \int_{\R^{2d-1}} \Big( \prod_{b=0}^1  \Big| \int_{\R} \Big( \prod_{j\in Q_0, j_d=b} F_j(x^j) \Big) \psi_b(x_{d}^b-p) dx_d^b \Big|  \Big)\] 
\[\times \Big( \prod_{l=1}^{d-1} \rho_l(x_l^1-x_l^0)\Big) (dx^0_l)_{l=1}^{d-1}(dx^1_l)_{l=1}^{d-1} dp.\]
Applying the Cauchy-Schwarz inequality in all $x^0_l,x^1_l$, $l\ne d$,  and $p$ 
bounds this by
\[  \prod_{b=0}^1 \Big( \int_{\R^{2d-1}}   \Big| \int_{\R} \Big( \prod_{j\in Q_0, j_d=b} F_j(x^j) \Big) \psi_b(x_{d}^b-p) dx_d^b \Big|^2 \] 
\[\times \Big( \prod_{l=1}^{d-1} \rho_l(x_l^1-x_l^0)\Big) (dx^0_l)_{l=1}^{d-1}(dx^1_l)_{l=1}^{d-1} dp \Big)^{1/2}.\]
Expanding the square and writing the integral in $p$ again as convolution,  the term for $b\in\{0,1\}$  equals $\Lambda(K_b,\mathbf{F})^{1/2}$, in view of  the identity \eqref{form_equiv}.
\end{proof}

The next result makes use of Lemma~\ref{lem:cs} to obtain a bound on the averages when the kernel has a mean zero function.
\begin{lemma} 
    \label{lem:gen_decay}   
    Let $0\le L_1\le L_2\le d$ be integers  and assume that $L_1< d$.  
    For $1\le l \le d$, let  $\varphi_l:\R\to \R$,  be a Schwartz function. For $L_1+1\le l \le d$, let $\psi_{0,l},\, \psi_{1,l}:\R\to \R$ be Schwartz functions with integral zero, and assume that 
    \begin{equation}
        \label{conv_assumption}
        \varphi_l=\psi_{0,l}*\psi_{1,l}.
    \end{equation}
 There exists a constant $C>0$ such that for any tuple of the negative integers $(m_l)_{l=L_1+1}^d$,     $t\in [1,2]$,    positive integer $I$, increasing sequence of the integers $(k_i)_{i=0}^I$,  and $\mathbf{f}$ as in \eqref{char_var}, the following holds. 
If 
    \[\Phi(s) = \Big( \prod_{l=1}^{L_1} \varphi_l(s_l) \Big)\Big( \prod_{l=L_1+1}^{L_2} 
    \varphi_l(2^{-m_l}s_l)
    \Big)\Big( \prod_{l=L_2+1}^{d} 
    \varphi_l(2^{-m_l}(s_l-1))
    \Big) \]
    and $m'=\min(m_{L_1+1},\ldots, m_d)$, 
then 
        \[\sum_{i =0}^I \|A^\Phi_{2^{k_i}t}(\mathbf{f})\|_q^q \le C 2^{q m'}I^{1-\frac{q}{2}}.\]

\end{lemma}

\begin{proof}
Since 
\[\sum_{i =0}^I \|A^\Phi_{2^{k_i}t}(\mathbf{f})\|_q^q = 2^{qm'}\sum_{i =0}^I \|A^{\widetilde{\Phi}}_{2^{k_i}t}(\mathbf{f})\|_q^q\]
with
$\widetilde{\Phi}  =     2^{-m'}\Phi$,
it suffices to show
 \[\sum_{i =0}^I \|A^{\widetilde{\Phi}}_{2^{k_i}t}(\mathbf{f})\|_q^q \le C I^{1-\frac{q}{2}}.\]

Moreover, it suffices to show that there is a constant $C_0$ independent of the numbers $k_i$ and $m_l$,   such that for any $|\varepsilon_i|\le 1$ and a normalized tuple of Schwartz functions $\mathbf{F}$  as in \eqref{singint_assumption}, the estimate 
  \begin{equation}
       \label{longshort2_form}
       |\Lambda(K,\mathbf{F})| \le C_0
  \end{equation}
holds with  
\begin{equation*}
        K =   \sum_{i=0}^I \varepsilon_i   \widetilde{\Phi}_{(2^{k_i}t)}.
\end{equation*}
Once this is shown,   Lemma \ref{lem:kintchine} applied to the   
functions $(\max(C_0,1)^{-1}\widetilde{\Phi}_{(2^{k_i}t)})_{i=0}^I$ finishes the proof. 

Let $l'$ be the index   such that $m_{l'}=m'$. We  will apply Theorem \ref{thm:singint} and  we distinguish two cases depending on the position of the index $l'$.  

First we prove \eqref{longshort2_form} if  $L_1+1\le l' \le L_2$.  By Lemma \ref{lem:rescalingtp} applied with $a_l=t$ for   $1\le l \le d$ we may assume  $t=1$. 
By Lemma~\ref{lem:rescalingtp}  applied with  $a_l = 2^{m_l}$ for all $L_1+1\le l \le L_2$ and $a_l =1$ otherwise,  it also suffices to prove \eqref{longshort2_form} with $\widetilde{\Phi}$ replaced by 
    \[s\mapsto \Big( \prod_{l=L_1+1, l\ne l'}^{L_2} 2^{m_l} \Big ) \Big( \prod_{l=1}^{L_2} \varphi_l(s_l) \Big)
     \Big( \prod_{l=L_2+1}^{d} 
    \varphi_l(2^{-m_l}(s_l-1))
    \Big). \]
Estimating $\prod_{l=L_1+1, l\ne l'}^{L_2} 2^{m_l}\le 1,$
it suffices to prove \eqref{longshort2_form} with $K =   \sum_{i=0}^I \varepsilon_i   \kappa_i $, where
\begin{equation}\label{Ki}
\kappa_i(s) =  \Big(\prod_{l=1}^{L_2}(\varphi_l)_{(2^{k_i})}(s_l)\Big) \prod_{l=L_2+1}^{d}
2^{-k_i}\varphi_l(2^{-m_l-k_i}s_l - 2^{-m_l}).
\end{equation}

The functions $\varphi_l$ are Schwartz and therefore
\begin{equation}
    \label{phil_schwartz_bd}
    |\varphi_l(u)| \le c_0   \rho(u)  
\end{equation}
for a constant $c_0$ depending on the Schwartz seminorms of $\varphi_l$, where we denoted    
\[\rho(u) = (1+|u|^2)^{-100d}.\]

 For  
  each $L_2+1\le l \le d$ we estimate uniformly in $m_l\le 0$
\[
|\varphi_l(2^{-m_l}(u-1))| \le c_0 (1+2^{-2m_l}|u-1|^2)^{-100d}\]
\begin{equation}
    \label{phil_unif_bd}
    \le c_0 (1+|u-1|^2)^{-100d} \le 2^{200d}c_0  \rho(u) .
\end{equation}

We apply Lemma \ref{lem:cs} to each kernel $\kappa_i$ in \eqref{Ki}, and we apply it  with $\ell_0=l'$,  $\rho_l=2^{-200d}c_0^{-1}\rho_{(2^{k_i})}$ for   $l\ne l'$, and 
\[\sigma_0=(\psi_{0,l'})_{(2^{k_i})}, \quad \sigma_1 = (\psi_{1,l'})_{(2^{k_i})}. \]
 Then we apply an additional Cauchy-Schwarz inequality in the summation in $i$ and use $|\varepsilon_i|\le 1$. This   estimates
\begin{equation}
    \label{ls_case2_lambdaest}
    | \Lambda(K,\mathbf{F})  |\le   C \Lambda({K}_{0},\mathbf{F})^{1/2}\Lambda({K}_{1},\mathbf{F})^{1/2} ,
\end{equation}
where for $b\in \{0,1\}$, 
\begin{equation}
    \label{kernelK0}
    K_{b}(s) = \sum_{i=0}^I (\psi_{b,l'}^\star)_{(2^{k_i})}*(\psi_{b,l'})_{(2^{k_i})}(s_{l'})\prod_{l=1,l\neq l'}^d \rho_{(2^{k_i})}(s_l).
    \end{equation}
Note that there are no more parameters $m_l$,  and the Fourier transforms of the kernels $K_{b}$ are of the form $\eqref{K_apply}$ with $\varepsilon_k=1$ when $k=k_i$ and $\varepsilon_k=0$ otherwise. The desired bound now follows from Theorem  \ref{thm:singint}.

Now we prove \eqref{longshort2_form} when  $L_2+1\le l' \le d$. Here we will also apply  the Cauchy-Schwarz inequality, after which we will   estimate the terms with the functions $\varphi_l$, $l\ne l'$, as before, and additionally we will be able to remove the translation parameter at $l=l'$. This will be needed since  at the index $l'
$ we have an additional factor of $2^{-m_{l'}}$ (recall the definition of $\widetilde{\Phi}$).

By Lemma \ref{lem:rescalingtp} applied with $a_l=t$ for   $1\le l \le d$ we may assume  $t=1$.  
By Lemma~\ref{lem:rescalingtp}  applied with  $a_l = 2^{m_l}$ whenever $L_1+1\le l \le L_2$ or $l=l'$, and with   $a_l =1$ otherwise, it suffices to prove \eqref{longshort2_form} with $\widetilde{\Phi}$ replaced by 
    \[s\mapsto \Big( \prod_{l=L_1+1}^{L_2} 2^{m_l} \Big ) \Big( \prod_{l=1}^{L_2} \varphi_l(s_l) \Big)
     \Big( \prod_{l=L_2+1,l\ne l'}^{d} 
    \varphi_l(2^{-m_l}(s_l-1))
    \Big)\varphi_{l'}(s_l-2^{-m_{l'}}). \]
Estimating $\prod_{l=L_1+1}^{L_2} 2^{m_l}\le 1,$
it suffices to prove
\eqref{longshort2_form} with $K =   \sum_{i=0}^I \varepsilon_i   \kappa_i $,  where
\[\kappa_i(s) = \Big(\prod_{l=1}^{L_2}(\varphi_l)_{(2^{k_i})}(s_l)\Big) \Big( \prod_{l=L_2+1,l\neq l'}^{d}
2^{-k_i}\varphi_l(2^{-m_l-k_i}s_l - 2^{-m_l}) 
\Big) 
(T_{2^{-{m_{l'}}}}\varphi_{l'})_{(2^{k_i})}(s_{l'}).\]
Here and in the rest of this proof, for a function $\rho:\R\to \C$ we denote its translation by $a\in \R$  as  $T_a\rho(u) = \rho(u - a)$, so that $(T_{2^{-{m_{l'}}}}\varphi_{l'}) (u) = \varphi_{l'}(u-2^{-m_{l'}} ) $.

By \eqref{conv_assumption} we have 
\[T_{2^{-{m_{l'}}}}\varphi_{l'} = T_{2^{-{m_{l'}}}}\psi_{0,l'}*\psi_{1,l'} \]
We again use the estimates \eqref{phil_schwartz_bd} and \eqref{phil_unif_bd} for $l\neq l'$. Then we    apply  Lemma \ref{lem:cs} to each kernel $\kappa_i$ with $\ell_0=l'$,    $\rho_l(u) = 2^{-200d}c_0^{-1}\rho_{(2^{k_i})}$ for  $l\ne l'$, and 
\[\sigma_0=(T_{2^{-{m_{l'}}}}\psi_{0,l'})_{(2^{k_i})}, \quad \sigma_1 = (\psi_{1,l'})_{(2^{k_i})}. \]
 This gives an estimate \eqref{ls_case2_lambdaest} with 
\[K_0(s)=\sum_{i=0}^I (T_{2^{-{m_{l'}}}}\psi_{0,l'})^\star_{(2^{k_i})}*(T_{2^{-{m_{l'}}}}\psi_{0,l'})_{(2^{k_i})}(s_{l'})\prod_{l=1,l\neq l'}^d \rho_{(2^{k_i})}(s_l)\]
and $K_1$ as in \eqref{kernelK0} for $b=1$. 
Since \[(T_{2^{-{m_{l'}}}}\psi_{0,l'})^\star_{(2^{k_i})}*(T_{2^{-{m_{l'}}}}\psi_{0,l'})_{(2^{k_i})} = (\psi_{0,l'})^\star_{(2^{k_i})}*(\psi_{0,l'})_{(2^{k_i})},\]
this kernel in fact equals \eqref{kernelK0} for $b=0$. The desired bound  now follows from Theorem ~\ref{thm:singint}.
\end{proof}

As a corollary we obtain the long variation bound when $L_1<d$. 
 \begin{lemma}
 Let $0\le L_1\le L_2\le d$ be integers  and assume that $L_1< d$.   
 There exists a constant $C>0$ such that for any tuple of negative integers $(m_l)_{l=L_1+1}^d$,  positive integer $I$, increasing sequence of the integers $(k_i)_{i=0}^I$, and  $\mathbf{f}$ as in \eqref{char_var}, the following holds. 
If 
    \[\Phi(s) = \Big( \prod_{l=1}^{L_1} \phi(s_l) \Big)\Big( \prod_{l=L_1+1}^{L_2} \phi_{0,m_l}(s_l) \Big) \Big(\prod_{l=L_2+1}^{d} \phi_{1,m_l}(s_l)\Big)  \]
    and $m'=\min(m_{L_1+1},\ldots, m_d)$, 
then 
 \[\sum_{i =1}^I \|A^\Phi_{2^{k_i}}(\mathbf{f}) - A^\Phi_{2^{k_{i-1}}}(\mathbf{f})\|_q^q \le C 2^{qm'}I^{1-\frac{q}{2}}.\]
\end{lemma}
Note that $q/q'\le q$ and thus $2^{qm'}\le 2^{\frac{q}{q'}m'}$, so we obtain the desired long variation bound \eqref{eq:long}.

\begin{proof}
By the triangle inequality for the $L^q$ and $\ell^q$ norms, it suffices to show
 \[\sum_{i =0}^I \|A^\Phi_{2^{k_i}}(\mathbf{f})\|_q^q \le C 2^{qm'}I^{1-\frac{q}{2}}.\]
 Expanding out the definitions of $\phi_{0,m}$ and $\phi_{1,m}$ in \eqref{phi1m}, 
    this follows   from Lemma~\ref{lem:gen_decay} applied with $t=1$,   $\varphi_l=\phi$ for $1\le l \le L_1$, and $\varphi_l=\widetilde{\theta}=\widetilde{\vartheta}*\vartheta$ for $L_1+1\le l \le d$. 
\end{proof}

\subsection{Short variation: estimate \texorpdfstring{\eqref{eq:short}}{}}
The following lemma relies on the fundamental theorem of calculus in $t$ and will be used to bound the short variation.
 \begin{lemma}\label{lem:ftc}
Let $p> 1$. There exists a constant $C>0$ such that   for any continuously differentiable $a: [2^k,2^{k+1}] \to \C$, $k\in \mathbb{Z}$,  positive integer $n$, and any   increasing sequence of real numbers  $(t_\ell)_{\ell=0}^L$ in $[2^k,2^{k+1}]$,   
\begin{equation}
    \label{ftc2_ptwse}
  \sum_{\ell=0}^n  |a(t_\ell)-a(t_{\ell-1})|^p \le C  \Big(  \int_{1}^{2} |a(2^kt)|^{p}\frac{dt}{t}\Big)^{1/p'} \Big(  \int_{1}^{2} |2^kta'(2^kt)|^p   \frac{dt}{t} \Big)^{1/p} ,
\end{equation}
where   $\frac{1}{p} + \frac{1}{p'}=1$. Moreover,
    \begin{equation}
\sum_{\ell=0}^n|a(t_\ell)-a(t_{\ell-1})|^p \leq  \int_1^{2} |2^kta'(2^kt)|^{p}\, \frac{dt}{t}.
\label{short1}
\end{equation}
\end{lemma}

\begin{proof}
(1) Proof of  \eqref{short1}. 
 By the fundamental theorem of calculus in  $t$,   
\[
\big|a(t_\ell)-a(t_{\ell-1})\big|^p =\Big|\int_{t_{\ell-1}}^{t_\ell}ta'(t)\frac{dt}{t}\Big|^p\,  \]
\[\le (t_\ell-t_{\ell-1})^{p-1} \int_{t_{\ell-1}}^{t_\ell} \big|ta'(t)\big|^p\, \frac{dt}{t^p} \le  \int_{t_{\ell-1}}^{t_\ell} \big|ta'(t)\big|^{p}\, \frac{dt}{t},
\]
where the first inequality follows from Jensen's inequality, and the second one follows by using $t \geq t_{\ell-1}\ge 2^k$, and $t_\ell-t_{\ell-1}\le 2^{k}$.  

Summing this inequality  over $\ell$ and using the disjointness of $[t_{\ell-1},t_\ell)\subseteq [2^{k},2^{k+1}]$, 
\[\sum_{\ell=0}^L|a(t_\ell)-a(t_{\ell-1})|^p \leq  \int_{2^k}^{2^{k+1}} |ta'(t)|^{p}\, \frac{dt}{t}.\]
A change of variables in $t$  now gives \eqref{short1}.  \\

(2) Proof of \eqref{ftc2_ptwse}. 
We will prove this under the assumption that $a$ is absolutely continuous,  real,  and non-negative. The desired inequality then follows by splitting $a$   into real and imaginary,   negative and positive parts, and use the triangle inequality for the $L^p$ and $\ell^p$ norms.

We will use the inequality 
   \begin{equation}
       \label{inequality}
   |a-b|^p \le |a^p-b^p|
   \end{equation}
   valid for any  $a,b>0$, whose proof we postpone to the end of the section. 
   
Using  \eqref{inequality} we 
  bound
\[ \big|a(t_\ell)-a(t_{\ell-1})\big|^{p} \le  \big|a(t_\ell)^p-a(t_{\ell-1})^p\big| . \]
Applying the fundamental theorem of calculus in variable $t$ gives  
\[\big|a(t_\ell)^p-a(t_{\ell-1})^p\big| = 
  \Big|\int_{t_{\ell-1}}^{t_\ell}-t\partial_t (a(t)^p ) \frac{dt}{t}\Big| =  p  \Big|\int_{t_{\ell-1}}^{t_\ell} a(t)^{p-1} ta'(t)
  \frac{dt}{t}\Big|  \]
    \[\le p \int_{t_{\ell-1}}^{t_\ell}|a(t)|^{p-1}|ta'(t)| \frac{dt}{t} . \]
    By  H\"older's inequality for the exponents $(p',p)$ and using that $(p-1)p' =p$, this is
    \[\le  p \Big(  \int_{t_{\ell-1}}^{t_\ell} |a(t)|^{p}\frac{dt}{t}\Big)^{1/p'} \Big(  \int_{t_{\ell-1}}^{t_\ell} |ta'(t)|^p   \frac{dt}{t} \Big)^{1/p}. \]
    Summing over $\ell$,  applying H\"older's inequality for the exponents $(p',p)$  in the summation, and using the disjointness of $[t_{\ell-1},t_\ell)$, we obtain
    \[\sum_{\ell=0}^n  |a(t_\ell)-a(t_{\ell-1})|^p \le C  \Big(  \int_{2^k}^{2^{k+1}} |a(t)|^{p}\frac{dt}{t}\Big)^{1/p'} \Big(  \int_{2^k}^{2^{k+1}} |ta'(t)|^p   \frac{dt}{t} \Big)^{1/p}. \]
   Changing variables in $t$ in each of the factors on the right-hand side   gives \eqref{ftc2_ptwse}.

It remains to prove \eqref{inequality}. 
Without loss of generality we may assume $a>b$. By homogeneity we may also assume $b=1$. Then it suffices to show that for any $a>1$,
 \[(a-1)^p \le a^p-1.\]

To see this, let $f(a)=(a-1)^p-a^p+1$. Since $f'(a)=p((a-1)^{p-1}-a^{p-1})\le 0$ the function is decreasing, so from $f(1)=0$ we get $f(a)\le f(1)=0$. This implies $(a-1)^p\le a^p-1$ for all $a > 1$ and $p>1$.   \end{proof}

 The following lemma gives the key estimate for the short variation bound. 
\begin{lemma}  
\label{lem:avg_B} 
Let $0\le L_1\le L_2\le d$ be integers. There is a constant $C>0$ such that for   
any tuple of the negative integers $(m_l)_{l=L_1+1}^d$,     $t\in [1,2]$, positive integer    $I$,   increasing sequence of the integers      $(k_i)_{i=0}^I$,  and     $\mathbf{f}$ as in \eqref{char_var}, the following holds. Let 
\begin{equation}
    \label{psi_gen}
    \Psi(s) =  \sum_{l_0=1}^d \psi_{l_0,m_{l_0}}(s_{l_0})\prod_{\substack{l=1,l\neq l_0}}^{d}\varphi_{l,m_l}(s_l), 
\end{equation}
where 
 $\varphi_{l,m_l} =\phi$ for $1\le l \le L_1$, $\varphi_{l,m_l} = \phi_{0,m_l}$ for $L_1+1\le l\le L_2$,    $\varphi_{l,m_l} = -\phi_{1,m_l}$ for $L_2+1\le l \le d$, and $\psi_{l,m_l}(u)= (u\varphi_{l,m_l}(u))'$ for all $1\le l \le d$. 
Then, 
          \[\sum_{i=0}^I  \|A_{2^{k_i}t}^\Psi(\mathbf{f})\|_q^{q} \le CI^{1-q/2}.\]
\end{lemma}

 \begin{proof}  We write
      $\Psi =  \sum_{l_0=1}^d\Psi_{l_0}$
      with
      \[\Psi_{l_0} (s) = \psi_{l_0,m_{l_0}}(s_{l_0})\Theta((s_l)_{l\ne l_0}), \]
where
\[\Theta((s_l)_{l\ne l_0}) = \prod_{\substack{l=1,l\neq l_0}}^{d}\varphi_{l,m_l}(s_l) \]
      \[=         \Big( \prod_{l=1,l\ne l_0}^{L_1} \phi(s_l) \Big)\Big( \prod_{l=L_1+1,l\ne l_0}^{L_2} \widetilde{\theta}(2^{-m_l}s_l)
      \Big) \Big( \prod_{l=L_2+1,l\ne l_0}^{d} \widetilde{\theta}(2^{-m_l}(s_l-1)) \Big )
      .\]
Here we have also expanded out the definitions  of the functions $\varphi_{l,m_{l}}$ and $\psi_{l,m_{l}}$  in \eqref{psi_gen}. 
By the triangle inequality,  it suffices to show 
\begin{equation}
    \label{est_triangle}
    \sum_{i=0}^I  \|A_{2^{k_i}t}^{\Psi_{l_0}}(\mathbf{F})\|_q^{q} \le CI^{1-q/2}.
\end{equation}

If $L_1=d$, then the last two products in $\Theta$ are empty and thus
\[\Psi_{l_0}(s) = (X\phi)'(s_{l_0})\prod_{l=1,l\ne l_0}^{d} \phi(s_l) , \]
where for a function $\rho$ on $\R$ we denote $(X\rho)(u) = u\rho(u)$, so that 
\[(X\rho)'(u) = (u\rho(u))' = \rho(u) +u\rho'(u).\]
 It suffices to show that there is a constant $C_0$ independent of the numbers $k_i, \,m_l,\, t$,  such that for any   $|\varepsilon_i|\le 1$  and a normalized tuple of Schwartz functions $\mathbf{F}$ as in \eqref{singint_assumption}, 
  \begin{equation}\label{shortvar_lambda}
       |\Lambda(K,\mathbf{F})| \le C_0, 
  \end{equation}
where  \[K =   \sum_{i=0}^I \varepsilon_i   (\Psi_{\ell_0})_{(2^{k_i}t)} .\]
An application of  Lemma \ref{lem:kintchine} to the sequence of functions $(\max(C_0,1)^{-1}(\Psi_{\ell_0})_{(2^{k_i}t)})_{i=0}^I$, together with homogeneity of the $L^q$ norm, then completes the proof. 

By Lemma \ref{lem:rescalingtp} it suffices to assume $t=1$.    Since $(X\phi)'$ is a Schwartz function with  integral zero, the kernel $K$ satisfies \eqref{symest} and the   bound \eqref{shortvar_lambda} follows from  Theorem~\ref{thm:singint}. 
 
So from  now on we assume $L_1 <d$. 
If $1\le l_0\le L_1$, then  again    
\[\psi_{l_0,m_{l_0}}(u) =  (X\phi)'(u).\]
Expanding out the definitions of $\phi_{0,{m_l}}$ and $\phi_{1,{m_l}}$, 
    the desired estimate \eqref{est_triangle} follows   from Lemma \ref{lem:gen_decay} applied with $\varphi_{l_0} = (X\phi)'$,  $\varphi_l=\phi$ for $l\ne l_0$, $1\le l \le L_1$, $\varphi_l=\widetilde{\theta}$ for $L_1+1\le l \le d$,  and $\psi_{l,0}=\widetilde{\vartheta},\, \psi_{l,1} = \vartheta$. 
    Note that Lemma \ref{lem:gen_decay}  also gives a factor of $2^{qm'}$, which we crudely estimate by one.

If $L_1+1\le l_0 \le L_2$, we compute with $m=m_{l_0}$
\[\psi_{l_0,m}(u)=\big(X\phi_{0,m}\big)'(u) =\phi_{0,m}(u) + u (\phi_{0,m})'(u)\]
\[=\widetilde{\theta}(2^{-m}u) +  2^{-m}u\widetilde{\theta}'(2^{-m}u))
 = (X\widetilde{\theta})'(2^{-m}u).\]
Further, using that  $\widetilde{\theta}$ is a primitive of $\theta=\vartheta*\vartheta$,  
\[(X\widetilde{\theta})'(u) =  \widetilde{\theta} + u\theta(u)\]
\begin{equation}
\label{xtheta_conv}
    =\widetilde{\theta}(u) + (X\vartheta)*\vartheta(u)+ \vartheta*(X\vartheta)(u).
\end{equation}
Note that   $X\vartheta$ is in frequency supported away from the origin.

We write $\Psi_{l_0} = \sum_{n=1}^3\Psi_{l_0,n}$, where 
\[\Psi_{l_0,1}(s) =  \widetilde{\theta} (2^{-m_{l_0}}s_{l_0})\Theta((s_l)_{l\ne l_0}),\]
\[\Psi_{l_0,2}(s) = ({X\vartheta}*\vartheta)(2^{-m_{l_0}}s_{l_0})\Theta((s_l)_{l\ne l_0}),\]
\[\Psi_{l_0,3}(s) =  ({\vartheta}*X\vartheta) (2^{-m_{l_0}}s_{l_0})\Theta((s_l)_{l\ne l_0}).\]
By the triangle inequality it suffices to prove 
    \[\sum_{i=0}^I  \|A_{2^{k_i}t}^{\Psi_{l_0,n}}(\mathbf{f})\|_q^{q} \le CI^{1-q/2}\] 
for $1\le n \le 3$.

These estimates  follow  by applying Lemma \ref{lem:gen_decay}   three times. If $n=1$ we apply  Lemma~\ref{lem:gen_decay}  with     $\varphi_l=\phi$ for  $1\le l \le L_1$, $\varphi_l=\widetilde{\theta}=\widetilde{\vartheta}*\vartheta$ for  $L_1+1\le l \le d$.
If $n=2$ we apply it with $\varphi_{l_0} =(X\vartheta)*\vartheta$, and if $n=3$ with $\varphi_{l_0} =\vartheta*(X\vartheta)$, while the functions
$\varphi_l$ for $l\ne l_0$  remain the same as for $n=1$. 
  Note that Lemma \ref{lem:gen_decay} again gives  a factor of $2^{qm'}$, which we crudely estimate by one.

It remains to treat the case $L_2+1\le l_0 \le d$.  With $m=m_{l_0}$  we compute
\[\psi_{l_0,m}(u)=-\big(X\phi_{1,m}\big)'(u) = - \phi_{1,m}(u) - u (\phi_{1,m})'(u)\]
 \[=\widetilde{\theta}(2^{-m}(u-1)) + 2^{-m}u\widetilde{\theta}'(2^{-m}(u-1))\]
\[=\widetilde{\theta}(2^{-m}(u-1)) + {(X\theta)}(2^{-m}(u-1)) + 2^{-m}{\theta}(2^{-m}(u-1))\]
\[=(X\widetilde{\theta})'(2^{-m}(u-1))+2^{-m}{\theta}(2^{-m}(u-1))\]
We use \eqref{xtheta_conv} and write $\Psi_{l_0} = \sum_{n=1}^4\Psi_{l_0,n}$, where 
\[\Psi_{l_0,1}(s) =  \widetilde{\theta} (2^{-m_{l_0}}(s_{l_0}-1))\, \Theta((s_l)_{l\ne l_0})\]
\[\Psi_{l_0,2}(s) =  ({X\vartheta}*\vartheta) (2^{-m_{l_0}}(s_{l_0}-1))\, \Theta((s_l)_{l\ne l_0})\]
\[\Psi_{l_0,3}(s) =  ({\vartheta}*X\vartheta) (2^{-m_{l_0}}(s_{l_0}-1))\, \Theta((s_l)_{l\ne l_0})\]
\[\Psi_{l_0,4}(s) = 2^{-m_{l_0}}\theta(2^{-m_{l_0}}(s_{l_0}-1))\, \Theta((s_l)_{l\ne l_0}).\]

We now apply  Lemma \ref{lem:gen_decay} four times. 
If $n=1$ we apply it with     $\varphi_l=\phi$ for  $1\le l \le L_1$, $\varphi_l=\widetilde{\theta}=\widetilde{\vartheta}*\vartheta$ for  $L_1+1\le l \le d$.
If $n=2$ we apply it with $\varphi_{l_0} =(X\vartheta)*\vartheta$,   if $n=3$ with $\varphi_{l_0} =\vartheta*(X\vartheta)$, and if $n=4$ we apply it with $\varphi_{l_0} =\theta = \vartheta*\vartheta$. 
If $n=2,3,4$, the  functions
$\varphi_l$ for $l\ne l_0$  remain the same  as for $n=1$. For $1\le n\le 3$  we obtain 
    \[\sum_{i=0}^I  \|A_{2^{k_i}t}^{\Psi_{l_0,n}}(\mathbf{F})\|_q^{q} \le C2^{qm'}I^{1-q/2} \le  CI^{1-q/2},\] 
 while  for $n=4$,
   \[\sum_{i=0}^I  \|A_{2^{k_i}t}^{\Psi_{l_0,4}}(\mathbf{F})\|_q^{q} \le C2^{q(m'-m_{l_0})}I^{1-q/2} \le  CI^{1-q/2}\] 
   since $m'-m_{l_0}\le 0$. 
The final bound for the averages associated with $\Psi_{l_0}$ now follows by the triangle inequality. 
\end{proof}

 Now we are ready to complete the short variation estimate. First we consider $L_1=d$.  
\begin{lemma} 
 \label{lem:short_var} 
 Let $\Phi = \prod_{l=1}^d \phi$.
There exists  a constant $C>0$ such that for  any     positive integer  $I$, increasing sequence of the  integers $(k_i)_{i=0}^I$,  and $\mathbf{f}$ as in \eqref{char_var},  
\begin{equation*}
    \sum_{i=0}^I \sup_{\substack{n\in \Z_+\\ 2^{k_i} < t_0<\ldots<t_n\le 2^{k_{i}+1}}} \sum_{\ell=1}^n \|A^\Phi_{t_\ell}(\mathbf{f})-A^\Phi_{t_{\ell-1}}(\mathbf{f})\|_q^q\leq C I^{1-\frac{q}{2}}.
\end{equation*}
 \end{lemma}

\begin{proof}  We have  $-t\partial_t(\Phi_t)=\Psi_t$, where $\Psi$ is as in \eqref{psi_gen} with $L_1=d$, i.e.  
\[    \Psi(s) =  \sum_{l_0=1}^d (X\phi)'(s_{l_0})\prod_{\substack{l=1,l\neq l_0}}^{d}\phi(s_l),  \]
This implies
\begin{equation*}
    -t\partial_t(A_t^\Phi(\mathbf{f}))=A_t^{\Psi}(\mathbf{f}).
\end{equation*}
Part \eqref{short1} of  Lemma \ref{lem:ftc} applied with $p=q$ and $a(t) = A_t^\Phi(\mathbf{f})(x)$ for each  fixed  $x$ gives
  \[\sum_{\ell=1}^n|A_{t_\ell}^\Phi(\mathbf{f})(x)-A_{t_{\ell-1}}^\Phi(\mathbf{f})(x)|^q \leq  \int_1^{2} |A_{2^{k_i}t}^{\Psi}(\mathbf{f})|^{q}\, \frac{dt}{t}\]
  for any $ t_0<\ldots <t_n \in (2^{k_i},2^{k_{i}+1}]$. 
 Integrating in $x$, taking the supremum,  summing in $i$, and using $t\ge 1$,  
gives 
\[   \sum_{i=0}^I \sup_{\substack{n\in \Z_+\\ 2^{k_i} < t_0<\ldots<t_n\le 2^{k_{i}+1}}} \sum_{\ell=1}^n \|A^\Phi_{t_\ell}(\mathbf{f})-A^\Phi_{t_{\ell-1}}(\mathbf{f})\|_q^q \le  \int_1^2 \sum_{i=0}^I \|A_{2^{k_i}t}^{\Psi}(\mathbf{f}) \|_q^q \, dt. \]
By Lemma \ref{lem:avg_B} there is a constant $C$  such that   the last display is no greater than 
\[\int_1^2 C I^{1-\frac{q}{2}} dt = C I^{1-\frac{q}{2}},\]
as desired.  
\end{proof}

Finally we bound $L_1<d$. 

\begin{lemma}
 \label{lem:short_var2} 
Let $0\le L_1\le L_2\le d$ be integers and assume $L_1<d$. There is a constant $C>0$ such that for   
any tuple of the negative integers $(m_l)_{l=L_1+1}^d$,      positive integer  $I$, increasing sequence of the  integers $(k_i)_{i=0}^I$,   and  $\mathbf{f}$ as in \eqref{char_var}, the following holds. Let
    \[\Phi(s) = \Big( \prod_{l=1}^{L_1} \phi(s_l) \Big)\Big( \prod_{l=L_1+1}^{L_2} \phi_{0,m_l}(s_l) \Big) \Big(\prod_{l=L_2+1}^{d} \phi_{1,m_l}(s_l)\Big)  \]
    and $m'=\min(m_{L_1+1},\ldots, m_d)$. Then 
\[   \sum_{i=0}^I  \sup_{\substack{n\in \Z_+\\ 2^{k_i} < t_0<\ldots<t_n\le 2^{k_{i}+1}}}  \sum_{\ell=1}^n \|A^\Phi_{t_\ell}(\mathbf{f})-A^\Phi_{t_{\ell-1}}(\mathbf{f})\|_q^q\le C2^{\frac{q}{q'}m'}I^{1-\frac{q}{2}} .\]
\end{lemma}

\begin{proof} We again use
$-t\partial_t(A_t^\Phi(\mathbf{f}))=A_t^{\Psi}(\mathbf{f})$, 
where $\Psi$ is now as in \eqref{psi_gen} with $L_1<d$. 
Applying \eqref{ftc2_ptwse} of Lemma \ref{lem:ftc} with $p=q$ and $a(t)= A^\Phi_t(\mathbf{f})(x)$ for each fixed $x$ gives 
\[\sum_{\ell=1}^n |A_{t_\ell}^\Phi(\mathbf{f})(x)-A_{t_{\ell-1}}^\Phi(\mathbf{f})(x)|^{q}   \le C  \Big( \int_{1}^2 |A_{2^{k_i}t}^\Phi(\mathbf{f})(x)|^{q}\frac{dt}{t}\Big)^{1/q'} \Big( \int_{1}^2 |A_{2^{k_i}t}^\Psi(\mathbf{f})|^q   \frac{dt}{t} \Big)^{1/q} . \]
 for any $ t_0<\ldots <t_n \in (2^{k_i},2^{k_{i}+1}]$. 
Integrating in $x$, taking the supremum,    summing in $i$,  and applying  H\"older's inequality in $x$ and $i$, we estimate 
\[   \sum_{i=0}^I  \sup_{\substack{n\in \Z_+\\ 2^{k_i} < t_0<\ldots<t_n\le 2^{k_{i}+1}}}  \sum_{\ell=1}^n \|A^\Phi_{t_\ell}(\mathbf{f})-A^\Phi_{t_{\ell-1}}(\mathbf{f})\|_q^q\]
\[\le C  \Big( \int_{1}^2   \sum_{i=0}^I  \|A_{2^{k_i}t}^\Phi(\mathbf{f})\|_q^{q}\,{dt}  \Big)^{1/q'} \Big( \int_{1}^2  \sum_{i=0}^I  \|A_{2^{k_i}t}^\Psi(\mathbf{f})\|_q^{q}\,{dt}  \Big)^{1/q} . \]
Using Lemma \ref{lem:gen_decay}   with $\varphi_l=\phi$ for $1\le l \le L_1$ and $\varphi_l=\widetilde{\theta} = \widetilde{\vartheta}*\vartheta$ for $L_1+1\le l \le d$,   
we estimate the integrand in the first factor, while using 
Lemma \ref{lem:avg_B} we estimate the integrand in the second factor. This bounds the last display as  
\[\le C  \Big( \int_{1}^2   2^{q m'}I^{1-\frac{q}{2}} {dt}  \Big)^{1/q'} \Big( \int_{1}^2  I^{1-q/2} {dt}  \Big)^{1/q} = C2^{\frac{q}{q'}m'}I^{1-\frac{q}{2}},   \] 
where we also integrated in $t$ for the final bound. 
\end{proof}

We note that, since Lemma \ref{lem:ftc} was applied pointwise in the proofs of Lemmas \ref{lem:short_var} and \ref{lem:short_var2}, these proofs  actually yield stronger pointwise estimates for short variation.

\section{Ergodic averages: Proof of Theorem \ref{thm:ergodic} using Theorem \ref{thm:char_q}}
\label{sec:transitionergodic}
First we show that  Theorem \ref{thm:char_q} implies that for any
$\varrho>2$, there exists a constant $C>0$ such that for any  positive integer $I$,      increasing sequence of positive real numbers $(t_i)_{i=0}^I$, and  
  $\mathbf{f}=(f_j)_{j\in Q} \in L^{2^d}(\R^d)$  with $\|f_{j}\|_{2^d}=1$,
\begin{equation}
\sum_{i=1}^I \|A_{t_i}(\mathbf{f})-A_{t_{i-1}}(\mathbf{f})\|_q^\varrho \le C.
\label{eq_char_rho}
\end{equation}
If $d=1$, then $q=2$ and Theorem \ref{thm:char_q} also gives the bound with $\varrho=2$.  

To see \eqref{eq_char_rho}, take an integer $I>0$ and     an increasing sequence of positive real numbers $(t_i)_{i=0}^I$.  For $k>0$ we define 
\[\mathcal{I}_k = \{1\le i \le I: 2^{-k}< \|A_{t_i}(\mathbf{f})-A_{t_{i-1}}(\mathbf{f}) \|_{q}\le  2^{-k+1}\}.\]
Applying Theorem \ref{thm:char_q} and using the lower bound on the difference of the averages gives
\begin{equation*}
2^{-kq}|\mathcal{I}_k| \le \sum_{i\in \mathcal{I}_k}\|A_{t_i}(\mathbf{f})-A_{t_{i-1}}(\mathbf{f}) \|_{q}^q \le C |\mathcal{I}_k|^{1-q/2},
\end{equation*}
which  implies $|\mathcal{I}_k| \le C 2^{2k}$.
With this  we then estimate for any $\varrho>2$
\[\sum_{i=1}^I \|A_{t_i}(\mathbf{f})-A_{t_{i-1}}(\mathbf{f})\|_{q}^\varrho = \sum_{k\ge 0}\sum_{i\in \mathcal{I}_k} \|A_{t_i}(\mathbf{f})-A_{t_{i-1}}(\mathbf{f})\|_{q}^\varrho \]

\[\le \sum_{k\ge 0} 2^{\varrho(-k+1)} |\mathcal{I}_k| \le C \sum_{k\ge 0}2^{(2-\varrho)k} \le C.\]

The transition of \eqref{eq_char_rho} to ergodic averages is  standard and is described for instance in \cite{T08}. Here we are using a more straightforward approach, as in \cite{DKST16}, and we give only a sketch of the argument.

First  we transfer the estimate \eqref{eq_char_rho} to $\Z^d$. 
For   functions $\tilde{f}_j\colon\Z^d\to\C$, $j\in Q$, normalized so that $\|\tilde{f}_j\|_{\ell^{2^d}(\Z^d)} = 1$,  we define the averages
\begin{equation*}
\tilde{A}_n(\mathbf{\tilde{f}})(\mathbf{k}) = \frac{1}{n^d}\sum_{i_1,\dots, i_d=0}^{n-1} \prod_{j\in Q}\tilde{f}_{j}(\mathbf{k}+j\odot\mathbf{i}), 
\end{equation*}
where $n\in\N$, $\mathbf{i}=(i_1,\ldots, i_d)$,  and $\mathbf{k}\in\Z^d$.

If for each $j\in Q$,  
$f_j:\R^d\to\R$ are given   by 
\begin{equation*}
f_j(x)=\tilde{f}_j(\lfloor x_1\rfloor,\dots, \lfloor x_d\rfloor)=\sum_{\mathbf{i}\in\Z^d}\tilde{f}_j(\mathbf{i})\mathbf{1}_{[i_1,i_1+1)}(x_1)\cdots\mathbf{1}_{[i_d,i_d+1)}(x_d),
\end{equation*}
then 
\begin{equation*}
\|f_{j}\|_{2^d}=\|\tilde{f}_{j}\|_{\ell^{2^d}(\Z^d)}=1,
\end{equation*}
and \[{A}_{n}(\mathbf{{f}})(\mathbf{k})= \tilde{A}_{n}(\mathbf{\tilde{f}})(\mathbf{k}).\]
Moreover,  for ${\alpha}=(\alpha_1,\ldots, \alpha_d)\in[0,1)^d$, 
\[{A}_{n}(\mathbf{{f}})(\mathbf{k}+{\alpha}) = \frac{1}{n^d}\sum_{i_1,\dots, i_d=0}^{n-1} \Big( \prod_{s=1}^da_{i_s} \Big ) \prod_{j\in Q}\tilde{f}_{j}(\mathbf{k}+j\odot\mathbf{i}), \]
where  for $s\in \{1,\ldots, d\}$, $a_{i_s} = 1-\alpha_s$ if $i_s=0$, $a_{i_s}=1$ if $1\le i_s \le n-2$, and $a_{i_s}=\alpha_s$ if $i_s=n-1$. 
Consequently, 
\begin{equation*}
\big | \|A_{n_i}(\mathbf{f})-A_{n_{i-1}}(\mathbf{f})\|_{q}
- \|\tilde{A}_{n_i}(\mathbf{\tilde{f}})-\tilde{A}_{n_{i-1}}(\mathbf{\tilde{f}})\|_{\ell^{q}(\Z^d)} \big |
\le  2^{d+1}(n_{i-1})^{-1}.
\end{equation*}
This, together with the triangle inequality, the estimate \eqref{eq_char_rho}, 
and $\sum_{i=1}^I n_{i-1}^{-\varrho} \le C$, gives
\begin{equation}
\sum_{i=1}^{I}\|\tilde{A}_{n_i}(\mathbf{\tilde{f}})-\tilde{A}_{n_{i-1}}(\mathbf{\tilde{f}})\|_{q}^{\varrho} \le C.
\label{var_za}
\end{equation}

Now we transfer to the probability space $(X,\mathcal{F},\mu)$. First let $f_j \in L^{\infty}(X)$, $j\in Q$, and normalize them as 
$\|f_j\|_{L^\infty(X)} =1.$
Take a point $x \in X$ and a positive integer $N\ge n_I$. Define the functions $\tilde{F}_{j}^{x,N}:\Z^d\to\C$, $j\in Q$, along the forward trajectory of $x$ by
\[\tilde{F}_{j}^{x,N}(\mathbf{k}) = \left\{\begin{array}{cl} f_{j}(T_1^{k_1}\cdots T_d^{k_d}x) & \text{ if } k_1,\dots,k_d\in\Z, \, 0\le k_1,\dots,k_d \le 2N-1\\ 
0 & \text{ otherwise. } 
\end{array}\right.\]
Since the transformations are commuting,
\begin{equation*}
M_n(\mathbf{f})(T_1^{k_1}\cdots T_d^{k_d}x)=\tilde{A}_n(\mathbf{\tilde{F}}^{x,N})(\mathbf{k}),
\end{equation*}
for all integers $0\le i_1,\dots,i_d <N$ and $0<n\le N$, where again $\mathbf{\tilde{F}}^
{x,N}$ is the $(2^d-1)$-tuple consisting of functions $\tilde{F}_{{j}}^{x,N}$. Since the transformations are measure-preserving, 
\begin{equation*}
\|M_{n_i}(\mathbf{f})-M_{n_{i-1}}(\mathbf{f})\|_{L^{q}(X)}^{q}\le  {N^{-d}}\int_X  \|\tilde{A}_{n_i}(\mathbf{\tilde{F}}^{x,N})-\tilde{A}_{n_{i-1}}(\mathbf{\tilde{F}}^{x,N}) \|_{\ell^{q}(\Z^d)}^{q}\,d\mu(x),
\end{equation*}
and then by Jensen's inequality since $\varrho>q$
\begin{equation*}
 \|M_{n_i}(\mathbf{f})-M_{n_{i-1}}(\mathbf{f})\|_{L^{q}(X)}^{\varrho}\le  {N^{-\frac{d\varrho}{q}}}\int_X  \|\tilde{A}_{n_i}(\mathbf{\tilde{F}}^{x,N})-\tilde{A}_{n_{i-1}}(\mathbf{\tilde{F}}^{x,N})\|_{\ell^{q}(\Z^d)}^{\varrho}\,d\mu(x).
\end{equation*}
Similarly, 
\begin{equation*}
\|\tilde{F}_j^{x,N}\|_{\ell^{2^d}(\Z^d)}\le (2N)^{\frac{d}{2^d}}\|f_j\|_{{L}^{\infty}(X)} = (2N)^{\frac{d}{2^d}}.
\end{equation*}
Applying \eqref{var_za} to the functions $\tilde{F}_{j}^{x,N}$ we obtain
\begin{equation}
\label{est_ergodic}
\sum_{i=1}^{I} \|M_{n_i}(\mathbf{f})-M_{n_{i-1}}(\mathbf{f}) \|_{L^{q}(X)}^{\varrho} \le C N^{-\frac{d\varrho}{q}}\Big(N^{\frac{d\varrho}{2^d}}\Big)^{2^d-1} = C
\end{equation}
for any $n_0<n_1<\cdots<n_I$, where we used   $q=2^d/(2^d-1)$.

To complete the proof of Theorem \ref{thm:ergodic}, we use the monotonicity and log-convexity of the $L^p$ norms. In the case $p< q$ we can use the monotonicity of $L^p$ norms on a probability space to get
\begin{equation*}
 \|M_{n_i}(\mathbf{f})-M_{n_{i-1}}(\mathbf{f})\|_{L^{p}(X)}\le \|M_{n_i}(\mathbf{f})-M_{n_{i-1}}(\mathbf{f})\|_{L^{q}(X)}.
\end{equation*}
 For $p>q$, by  log-convexity of $L^p$ norms, 
\[
\|M_{n_i}(\mathbf{f})-M_{n_{i-1}}(\mathbf{f})\|_{L^{p}(X)} 
\le \Big(2\prod_{j\in Q}\|f_j\|_{L^{\infty}(X)}\Big)^{1-\frac{q}{p}}\|M_{n_i}(\mathbf{f})-M_{n_{i-1}}(\mathbf{f})\|_{L^{q}(X)}^{\frac{q}{p}}\]
\[=2 \|M_{n_i}(\mathbf{f})-M_{n_{i-1}}(\mathbf{f})\|_{L^{q}(X)}^{\frac{q}{p}} .
\]
Taking the $\varrho$-th power, summing in $i$,   using \eqref{est_ergodic}, and in the case $p>q$ also the condition $q\varrho/p>2$, we finish the proof of Theorem \ref{thm:ergodic}.

\end{document}